\documentclass[11pt]{amsart}
\usepackage{amsfonts}
\usepackage{amssymb, amscd}
\usepackage{amsthm}
\usepackage{amsmath}
\usepackage{latexsym}

\newtheorem{thm}{Theorem}[section] 
\newtheorem{cor}[thm]{Corollary}
\newtheorem{defn}[thm]{Definition}
\newtheorem{example}[thm]{Example}
\newtheorem{lemma}[thm]{Lemma}

\newtheorem{remark}[thm]{Remark}

\newtheorem{prob}[thm]{Problem}

\numberwithin{equation}{section}

\newcommand{\OO}{{\mathcal O}}

\newcommand{\VV}{{\mathcal V}}
\newcommand{\FC}{{\mathcal F}}
\newcommand{\MC}{{\mathcal M}}
\newcommand{\IC}{{\mathcal IC}}

\newcommand{\HB}{\mathbb{H}}
\newcommand{\Z}{\mathbb{Z}}
\newcommand{\Lef}{\mathbb{L}}

\newcommand{\Q}{\mathbb{Q}}

\newcommand{\C}{\mathbb{C}}

\newcommand{\mh}{\mbox{MHM}}

\newcommand{\mhs}{\mbox{mHs}}

\begin{document}

\title{Nearby cycles and characteristic classes of singular spaces}

\author[J. Sch\"urmann ]{J\"org Sch\"urmann}
\address{J.  Sch\"urmann : Mathematische Institut,
          Universit\"at M\"unster,
          Einsteinstr. 62, 48149 M\"unster,
          Germany.}
\email {jschuerm@math.uni-muenster.de}

\begin{abstract} In this paper we give an introduction to our recent work on characteristic classes of complex hypersurfaces based on some talks given at conferences in Strasbourg, Oberwolfach and Kagoshima. We explain the relation between nearby cycles for constructible functions or sheaves as well as for (relative) Grothendieck groups of algebraic varieties and mixed Hodge modules, and the specialization of characteristic classes of singular spaces like 
the Chern-, Todd-, Hirzebruch- and motivic Chern-classes. As an application we get a description of the differences between the corresponding virtual and functorial characteristic classes of complex hypersurfaces in terms of vanishing cycles related to the singularities of the hypersurface.    
\end{abstract}

\maketitle

A natural problem in complex geometry  is the relation between invariants of a singular complex hypersurface $X$ (like  Euler characteristic and Hodge numbers)
and the geometry of the singularities of the hypersurface (like the local Milnor fibrations).
For the Euler characteristic this is for example a special case of the difference between the Fulton- and MacPherson- Chern classes of $X$, whose differences are the now well studied Milnor classes of $X$ (\cite{Alu0,BLSS,BLSS2,BSS,Max,PP,Sch,Sch2,Y}). Their degrees are related to Donaldson-Thomas invariants of the singular locus (\cite{Beh}).\\

A very powerful approach to this type of questions is by the theory of the nearby and vanishing cycle functors. For example a classical result of Verdier \cite{V} says that
the MacPherson Chern class transformation \cite{MP, Ken} commutes with specialization, which for constructible functions means the corresponding nearby cycles. Here we explain the corresponding result for our motivic Chern- and Hirzebruch class transformations 
as introduced in our joint work with J.-P. Brasselet and S. Yokura \cite{BSY},
i.e. they also commute with specialization defined in terms of nearby cycles.
Here one can work either in the motivic context with relative Grothendieck group of varieties
\cite{Bit, GLM}, or in the Hodge context with Grothendieck groups of M. Saito's mixed Hodge modules \cite{Sa0, Sa1}. The key underlying specialization  result \cite{Sc} is about the filtered de Rham complex of the underlying filtered $D$-module in terms of the Malgrange-Kashiwara $V$-filtration. But here we focus on the geometric motivations and applications as given in our joint work with S.E. Cappell, L. Maxim and J.L. Shaneson
\cite{CMSS}.\\

In this paper we work (for simplicity) in the complex algebraic context,
since this allows us to switch easily between an algebraic geometric language and
an underlying topological picture.
Many results are also true in the complex analytic or algebraic context over base field of characteristic zero. First we introduce the virtual characteristic classes and numbers of
hypersurfaces and local complete intersetions in smooth ambient manifolds. Next we recall some
of the theories of  functorial characteristic classes for singular spaces 
\cite{MP, BFM, CS1, BSY, Sch3}. Finally we explain the relation to nearby and vanishing cycles following our earlier results \cite{Sch, Sch2}
about different Chern classes for singular spaces.\\

{\em Acknowledgements:}
This paper is an extended version of some talks  given  
at conferences in Strasbourg, Oberwolfach and Kagoshima.
Here I would like to thank the organizers  for the invitation to these conferences.
I also would like to thank Sylvain Cappell, Laurentiu Maxim and Shoji Yokura for the discussions on 
our joint work related to the subject of this paper.

\tableofcontents

\section{Virtual classes of local complete intersections}

Recall that we are working  in the complex algebraic context.
A characteristic class $cl^*$ of (complex algebraic) vector bundles over $X$ is a map
$$cl^*: Vect(X) \to \HB^*(X) \otimes R$$
from the set $Vect(X)$ of isomorphism classes of complex algebraic vector bundles over $X$ to some  
cohomology theory $\HB^*(X) \otimes R$ with a coefficient ring $R$, which is compatible with pullbacks. Here we use as a cohomology theory 
$$\HB^*(X)= \begin{cases}
H^{2*}(X,\Z) &,\text{the usual cohomology in even degrees.}\\
CH^*(X) &,\text{the operational Chow cohomology of \cite{Fu}.}\\
K^0(X) &,\text{the Grothendieck group of vector bundles.}
\end{cases}$$
We also assume that $cl^*$ is {\em multiplicative}, i.e.
$$cl^*(V)=cl^*(V')\cup cl^*(V'')$$
for any short exact sequence 
$$0\to V' \to V\to V''\to 0$$
of vector bundles on $X$, with $\cup$ given by the cup- or tensor-product. 
Such a characteristic class $cl^*$ corresponds by the ``splitting principle''
to a unique formal power series $f(z)\in R[[z]]$ with
$cl^*(L)=f(c^1(L))$ for any line bundle $L$ on $X$. Here $c^1(L)\in \HB^1(X)$ is
the nilpotent first Chern class of $L$, which in the case $\HB^*(X)=K^0(X)$ is given by
$c^1(L):=1-[L^{\vee}]\in K^0(X)$ (with $(\cdot)^{\vee}$ the dual bundle).
Finally $cl^*$ should be {\em stable} in the sense that $f(0)\in R$ is a unit
so that $cl^*$ induces a functorial group homomorphism
$$cl^*: \left(K^0(X),\oplus\right)\to \left(\HB^*(X)\otimes R,\cup\right)\:.$$

Let us now switch to smooth manifolds, which will be an important intermediate step on the way to characteristic classes of singular spaces. For a  complex algebraic manifold $M$ its tangent bundle $TM$ is available and a characteristic class $cl^{*}(TM)$ of the tangent bundle $TM$ is called a {\em characteristic cohomology class} $cl^{*}(M)$ of the manifold $M$. We also use the notation 
$$cl_{*}(M):=cl^{*}(TM)\cap [M]\in \HB_{*}(M)\otimes R$$
for the corresponding {\em characteristic homology class} of the manifold $M$, with $[M]\in \HB_{*}(M)$ the fundamental class (or the class of the structure sheaf) in
$$\HB_{*}(M):=\begin{cases}
H^{BM}_{2*}(M) &,\text{the Borel-Moore homology in even degrees.}\\
CH_*(M) &,\text{the Chow group.}\\
G_0(M) &,\text{the Grothendieck group of coherent sheaves.}
\end{cases}$$

If $M$ is moreover compact, i.e. the constant map $k: M\to \{pt\}$ is proper,
one gets the corresponding {\em characteristic number}
\begin{equation}
\sharp(M):=k_*(cl_*(M))=:deg(cl_*(M))\in R\:.
\end{equation}

\begin{example}[Hirzebruch '54] The famous
Hirzebruch $\chi_y$-genus is the characteristic number whose associated characteristic class can be given in two versions (see \cite{H}):
\begin{enumerate}
\item The cohomological version, with $R=\Q[y]$, is given by the Hirzebruch class $cl^*=T^*_y$
corresponding to the normalized power series  
\[f(z):=Q_{y}(z):= \frac{z(1+y)}{1-e^{-z(1+y)}} -zy
\quad \in \Q[y][[z]] \:. \]
\item The $K$-theoretical version, with $R=\Z[y]$, is given by the {\em dual total Lambda-class} $cl^*=\Lambda_y^{\vee}$, with
\[\Lambda_y^{\vee}(\cdot):=\Lambda_y\left((\cdot)^{\vee}\right)=
\sum_{i\geq 0}\; \left[\Lambda^i \left((\cdot)^{\vee}\right)\right]\cdot y^i\]
corresponding to the unnormalized power series 
$$f(z)=1+y-yz\in \Z[y][[z]]\:.$$
\end{enumerate}
\end{example}

So the $\chi_y$-genus of the compact complex algebraic manifold $M$ is given by
\begin{align*}
\chi_{y}(M):= 
&\sum_{p\geq 0} \chi(M, \Lambda^{p}T^{*}M)\cdot y^{p}\\
=&\sum_{p\geq 0} \left( \sum_{i\geq 0}
(-1)^{i}dim_{C}H^{i}(M,\Lambda^{p}T^{*}M) \right)\cdot y^{p} \:,
\end{align*}
with $T^*M$ the algebraic cotangent bundle of $M$. The equality
\begin{equation} \label{eq:gHRR} 
\chi_{y}(M)= deg\left( T^{*}_{y}(TM) \cap [M]\right)
\quad \in \Q[y], \tag{gHRR}
\end{equation}
is (called) the {\em generalized Hirzebruch Riemann-Roch theorem\/}
\cite{H}. The corresponding power series $Q_{y}(z)$ (as above) specializes to
\begin{displaymath}
Q_{y}(z) = 
\begin{cases}
\:1+z &\text{for $y=-1$,}\\
\:\frac{z}{1-e^{-z}} &\text{for $y=0$,}\\
\:\frac{z}{\tanh z} &\text{for $y=1$.}
\end{cases} \end{displaymath}

Therefore the {\em Hirzebruch class\/} $T^{*}_{y}(TM)$ unifies the
following important (total) characteristic cohomology classes of $TM$:
\begin{equation}\label{imp-classes}
T^{*}_{y}(TM) = 
\begin{cases}
\:c^{*}(TM) &\text{the  {\em Chern class\/} for $y=-1$,}\\
\:td^{*}(TM) &\text{the  {\em Todd class\/} for $y=0$,}\\
\:L^{*}(TM) &\text{the  {\em Thom-Hirzebruch L-class\/} for $y=1$.}
\end{cases} \end{equation}

The gHHR-theorem  specializes  to
the calculation of the following important invariants:
\begin{equation}\label{imp-genera}
\begin{split}
\chi_{-1}(M) &= e(M) = deg\left( c^{*}(TM)\cap [M]\right) \quad \text{the {\em Euler characteristic\/},}\\ 
\chi_0(M) &= \chi(M) = deg\left( td^{*}(TM)\cap [M]\right) \quad\text{the {\em arithmetic genus\/},}\\ 
\chi_1(M) &= sign(M) = deg\left( L^{*}(TM)\cap [M]\right)  \quad\text{the {\em signature\/},}
\end{split} \end{equation}
which are, respectively, the {\em Poincar\'e-Hopf or Gauss-Bonnet theorem\/}, the
{\em Hirzebruch Riemann-Roch theorem\/} and the
{\em Hirzebruch signature theorem\/}.\\

If $X$ is a singular complex algebraic variety, then the algebraic tangent bundle
of $X$ doesn't exist so that a characteristic (co)homology class of $X$ can't be defined as before. But if $X$ can be realized as a {\em local complete intersection} inside a complex
algebraic manifold $M$, then a substitute for $TX$ is available. Indeed this just means
that the closed inclusion $i: X\to M$ is a {\em regular embedding} into the smooth algebraic manifold $M$, so that the {\em normal cone} $N_XM\to X$ is an algebraic vector bundle over $X$
(compare \cite{Fu}).  Then the {\it virtual tangent bundle} of $X$ 
\begin{equation}
T_{{\rm vir}}X:=[i^*TM-N_XM] \in K^0(X),
\end{equation} 
is {\em independent} of the embedding in $M$ (e.g., see \cite{Fu}[Ex.4.2.6]), so it is a well-defined element in the Grothendieck group of vector bundles on $X$. Of course 
$$T_{{\rm vir}}X=[TX] \in K^0(X)$$
in case $X$ is a smooth algebraic submanifold.\\

If $cl^*: K^0(X)\to \HB^*(X)\otimes R$ denotes a  characteristic cohomology class as before, then one can associate to $X$ an {\it intrinsic}  homology class (i.e., independent of the embedding $X\hookrightarrow M$) defined as:
\begin{equation}
cl_*^{\rm vir}(X):=cl^*(T_{{\rm vir}}X) \cap [X] \in \HB_*(X)\otimes R\:.
\end{equation} 
Here $[X]\in \HB_*(X)$ is again the {\it fundamental class} (or the class of the structure sheaf) of $X$ in  
$$\HB_{*}(X):=\begin{cases}
H^{BM}_{2*}(X) &,\text{the Borel-Moore homology in even degrees.}\\
CH_*(X) &,\text{the Chow group.}\\
G_0(X) &,\text{the Grothendieck group of coherent sheaves.}
\end{cases}$$
Here $\cap$ in the $K$-theoretical context comes from the tensor product with the coherent locally free sheaf of sections of the vector bundle. Moreover, for the class $cl^*=\Lambda_y^{\vee}$ one has to take $R:=\Z[y,(1+y)^{-1}]$ to make it a stable characteristic class defined on $K^0(X)$.\\

Let $i: X\to M$ be a regular embedding of (locally constant) codimension $r$ between possible  singular complex algebraic varieties. Using the famous {\em deformation to the normal cone},
one gets functorial {\em Gysin homomorphisms} (compare \cite{Fu, V0, V})
\begin{equation}
i^!: \HB_*(M)\to \HB_{*-r}(X)
\end{equation}
and 
\begin{equation}
i^!: G_0(M)\to G_0(X) \:.
\end{equation}
Note that $i$ is of finite tor-dimension, so that the last $i^!$ can also be described as
$$i^!=Li^*: G_0(M)\simeq K_0(D^b_{coh}(M))\to K_0(D^b_{coh}(X))\simeq G_0(X)$$
coming from the derived pullback $Li^*$ between the bounded derived categories with
coherent cohomology sheaves. If $M$ is also {\em smooth}, then one gets easily the following important relation between the virtual characteristic classes $cl^{vir}_*(X)$ of $X$ and the
Gysin homomorphisms:
\begin{equation}\label{gysin-vir}
i^!\left( cl_*(M)\right) = i^!\left( cl^*(TM)\cap [M]\right)=
cl^*(N_XM)\cap cl_*^{vir}(X) \:.
\end{equation}

From now on we assume that 
$$X=\{f=0\}=\{f_i=0\;|\;i=1,\dots,n\}$$ 
is a {\em global complete intersection} in the complex algebraic manifold $M$ coming from a cartesian diagram
\begin{equation} \begin{CD}
\{f=0\} @= X @> i >> M\\
@. @V f_0 VV @V f= V (f_1,\dots,f_n) V\\
@. \{0\} @> i_0 >> \C^n \:.
\end{CD}\end{equation}

Then $N_XM\simeq f^*\left(N_{\{0\}}\C^n\right)=X\times \C^n$ is a {\em trivial} vector bundle of rank $n$ on $X$ so that
\begin{equation}\label{factor-N}
cl^*\left(N_XM\right)=\begin{cases}
1 &\text{for $cl^*=T^*_y,\; c^*,\; td^*$ or $L^*$.}\\
(1+y)^n &\text{for $cl^*=\Lambda_y^{\vee}$.}
\end{cases}
\end{equation}

Assume now that $f$ is proper so that $X$ is compact. Since the Gysin homomorphisms $i^!$ commute with proper pushdown (compare \cite{Fu, V0, V}), one gets by the projection formula
\begin{align*}
\sharp^{vir}(X):=f_{0*}\left( cl^{vir}_*(X)\right)
&= f_{0*}\left( cl^*(N_XM)^{-1} \cap i^!cl_*(M)\right)\\
&=cl^*(N_{\{0\}}\C^n)^{-1} \;\cap\; i_0^!\left( f_*cl_*(M)\right)\:.
\end{align*}
Taking a (small) regular value $0\neq t \in \C^n$, in the same way
from the cartesian diagram
\begin{equation} \begin{CD}
\{f=0\} @= X @> i >> M @< i' << X_t @= \{f=t\} \\
@. @V f_0 VV @V f VV @VV f_t V\\
@. \{0\} @> i_0 >> \C^n @< i_t << \{t\}\
\end{CD}\end{equation}
for the ``nearby'' smooth submanifold $X_t=\{f=t\}$,  one gets the equality
$$\sharp(X_t):=f_{t*}\left( cl_*(X_t)\right)=
cl^*(N_{\{t\}}\C^n)^{-1} \;\cap\; i_t^!\left( f_*cl_*(M)\right)\:.$$
Note that the set of critical values of $f$ is a proper algebraic subset of
$\C^n$, as can be seen by ``generic smoothness'' or from an adapted stratification 
of the proper algebraic map $f$.
Now $N_{\{0\}}\C^n\simeq \C^n \simeq N_{\{t\}}\C^n$ and the {\em smooth pullback}
$\pi^*$ for the (vector bundle) projection $\pi: \C^n\to \{pt\}$ is an isomorphism
$$\pi^*: R=\HB_*(\{pt\})\otimes R \simeq \HB_{*+n}(\C^n)\otimes R$$
with inverse $i_0^!$ and $i_t^!$ (see \cite{Fu, V0, V}), so that
the ``virtual characteristic number''
\begin{equation}
\sharp^{vir}(X):=f_{0*}\left( cl^{vir}_*(X)\right)=\sharp(X_t) \in R
\end{equation}
is the corresponding characteristic number of a ``nearby'' smooth fiber $X_t$.

\section{Functorial characteristic classes of singular spaces}

For a more general singular complex algebraic variety $X$ its ``virtual tangent bundle'' is not available any longer, so characteristic classes for singular varieties have to be defined in a different way. For an introduction to this subject compare with our survey paper \cite{SY} (and see also \cite{Sch3, Yo}).
The theory of characteristic classes of vector bundles is a natural transformation of contravariant functorial theories. 
This {\em naturality} is an important guide for developing various theories of characteristic classes for singular varieties. Almost all known {\em characteristic classes\/} for singular spaces are formulated as {\em natural transformations\/} 
$$cl_{*}: A(X)\to \HB_{*}(X)\otimes R $$
of covariant functorial theories. Here  $A$ is a suitable theory (depending on the choice of $cl_{*}$),
which is covariant functorial for proper algebraic morphisms.\\

There is always a {\em distinguished element\/} $\mathrm{I}_{X}\in A(X)$ such that the corresponding {\em characteristic class of the singular space\/} $X$ is defined as
$$cl_{*}(X):=cl_{*}(\mathrm{I}_{X})\:.$$
Finally one has the {\em normalization\/}
$$cl_{*}(\mathrm{I}_{M})= cl^{*}(TM)\cap [M] \in \HB_{*}(M)\otimes R $$
for $M$ a smooth manifold, with $cl^{*}(TM)$ the corresponding characteristic cohomology class of $M$. 
This justifies the notation $cl_{*}$ for this homology class transformation, which should be seen as a homology class version of the following {\em characteristic number\/} of the singular space $X$:
$$\sharp (X):=cl_{*}(k_{*}\mathrm{I}_{X})= deg\left(cl_{*}(\mathrm{I}_{X})\right) \in 
\HB_{*}(\{pt\})\otimes R\simeq R \:,$$
with $k: X\to \{pt\}$ a constant map. Note that the {\em normalization\/} implies that for $M$ smooth:
$$\sharp (M)=deg\left(cl_{*}(M)\right)= deg\left(cl^{*}(TM) \cap [M]\right)$$
so that this is consistent with the notion of characteristic number of the smooth manifold $M$ as used before.\\

But only few characteristic numbers and classes have been extended in this way to singular
spaces. For example the three characteristic numbers (\ref{imp-genera}) and classes (\ref{imp-classes})
have been generalized to a singular complex
algebraic variety $X$ in the following way (where the characteristic numbers are only
defined for $X$ compact):

\begin{equation} \label{eq:y-1} \tag{$y=-1$}
 e(X) = deg\left( c_{*}(X)\right), \quad \text{with} \quad
c_{*}: F(X) \to H_{*}(X)
\end{equation}
the {\em Chern class transformation\/} of MacPherson \cite{MP, Ken}
from the abelian group $F(X)$
of complex algebraically constructible functions to homology,
where one can use the Chow group $CH_{*}(\cdot)$ or the Borel-Moore homology group
$H^{BM}_{2*}(\cdot,\Z)$ (in even degrees). Here $e(X)$ is the (topological) 
{\em Euler characteristic} of $X$, and the distinguished element $\mathrm{I}_{X}:=1_X\in F(X)$ is simply given by the characteristic function of $X$.
Then $c_{*}(X):=c_{*}(1_{X})$
agrees by \cite{BrS} via ``Alexander duality'' for compact $X$ embeddable into a complex manifold with
the  {\em Schwartz class\/} of $X$  as introduced before by
M.-H. Schwartz \cite{Schwa}.

\begin{equation} \label{eq:y0} \tag{$y=0$}
\chi(X) = deg\left(td_{*}(X)\right), \quad \text{with} \quad
td_{*}: G_{0}(X) \to H_{*}(X)\otimes \Q 
\end{equation}
the {\em Todd transformation\/} in the singular Riemann-Roch theorem of
Baum-Fulton-MacPherson \cite{BFM} (for Borel-Moore homology) or Fulton \cite{Fu}
(for Chow groups). Here $G_{0}(X)$ is the Grothendieck group
of coherent sheaves, with $\chi(X)$ the {\em arithmetic genus} (or holomorphic Euler characteristic) of $X$.
Then $td_{*}(X):=td_{*}([\OO_{X}])$,
with the distinguished element $\mathrm{I}_{X}:=[\OO_{X}]$ the class of the structure sheaf.\\

Finally for compact $X$ one also has
\begin{equation} \label{eq:y1} \tag{$y=1$}
sign(X) = deg\left( L_{*}(X)\right), \quad \text{with} \quad
L_{*}: \Omega(X) \to H_{2*}(X,\Q) 
\end{equation}
the {\em homology $L$-class transformation\/} of Cappell-Shaneson \cite{CS1}
as formulated in \cite{BSY}. 
Here $\Omega(X)$ is the abelian group
of cobordism classes of selfdual constructible complexes.
Then $L_{*}(X):=L_{*}([\IC_{X}])$ is the {\em homology L-class\/} of
Goresky-MacPherson \cite{GM1},
with the distinguished element $\mathrm{I}_{X}:=[\IC_{X}]$ the class of their intersection cohomology complex. So $sign(X)$ is the intersection cohomology signature of $X$.
For a {\em rational PL-homology manifold\/} $X$,
these $L$-classes are due to Thom \cite{Thom}.\\

So all these theories have the same formalism, but they are defined 
on completely different theories. Nevertheless, it is natural to ask
for another theory of characteristic homology classes of 
singular complex algebraic varieties, which unifies the above characteristic
homology class transformations. 
Of course in the smooth case, this is done by the 
{\em Hirzebruch class\/} $T^{*}_{y}(TM) \cap [M]$ of the tangent
bundle. An answer to this question was given in \cite{BSY} (together with some improvements in \cite{Sch3}). Using Saito's deep theory of {\em algebraic mixed Hodge modules} \cite{Sa0, Sa1},
we introduced in \cite{BSY} the {\em motivic Chern class transformations} as natural transformations (commuting with proper push down)
fitting into a commutative diagram:
\begin{equation*}
\begin{CD}
 G_0(X)[y] @>>> G_0(X)[y,y^{-1}] @= G_0(X)[y,y^{-1}] \\
@A mC_y AA @A mC_y AA @AA MHC_y A \\
K_0(var/X) @>>> {\mathcal M}(var/X) @> \chi_{Hdg}>> K_0(MHM(X)) \:.\\
\end{CD}
\end{equation*}

Here  $K_0(MHM(X))$ is the Grothendieck group of  algebraic mixed Hodge modules on $X$, and
$K_0(var/X)$ (resp. ${\mathcal M}(var/X):=K_0(var/X)[\Lef^{-1}]$)
is the (localization of the) relative  Grothen\-dieck group of complex algebraic varieties over $X$  
(with respect to the class of the affine line $\Lef$, compare e.g. \cite{Bit, GLM}).
The distinguished element is given by the constant Hodge module (complex) resp. by the class of the identity arrow
$$\mathrm{I}_{X}:=[\Q^H_X]\in K_0(MHM(X)) \quad \text{resp.} \quad
\mathrm{I}_{X}:=[id_X]\in K_0(var/X)\:,$$
and the canonical ``Hodge realization'' homomorphism $\chi_{Hdg}$ is given by
\begin{equation}
\chi_{Hdg}: K_0(var/X)\to K_0(MHM(X));\: [f: Y\to X]\mapsto [f_!\Q^H_Y] \:.
\end{equation}

The motivic Chern class transformations $mC_y, MHC_y$ capture information about 
the {\em filtered de Rham complex} of the filtered $D$-module underlying a mixed Hodge module. The corresponding characteristic class of the space $X$,
$$mC_y(X)=MHC_y(X) \in G_0(X)[y]\:,$$
can also be defined with the help of the (filtered) {\em Du Bois complex} of $X$ \cite{DB},
and satisfies for $M$ smooth the normalization condition
\begin{equation}\label{C-norm}
 mC_y(M)=MHC_y(M) = \Lambda_y^{\vee}(TM)\cap [M] \in G_0(M)[y]\:.
\end{equation}

The motivic Chern class transformations are
a $K$-theoretical refinement of the {\em Hirzebruch class transformations} $T_{y*}, MHT_{y*}$,
which can be defined by the (functorial) commutative diagram :
 \begin{equation*}
\begin{CD}
 {\mathcal M}(var/X) @> \chi_{Hdg}>> K_0(MHM(X)) @> MHC_y >>  G_0(X)[y,y^{-1}]\\
@V T_{y*} VV @V MHT_{y*} VV @VV td* V \\
H_*(X)\otimes \Q[y,y^{-1}] @>>>  H_*(X)\otimes\Q_{loc}   @< (1+y)^{-*}\cdot<<  H_*(X)\otimes \Q[y,y^{-1}] \:,
\end{CD}
\end{equation*}

with $td_*: G_0(X)\to H_*(X)\otimes\Q$ the {\em Todd class transformation} of Baum-Fulton-MacPherson \cite{BFM, Fu} and $(1+y)^{-*}\cdot$ the renormalization
given in degree $i$ by the multiplication
$$(1+y)^{-i}\cdot : H_i(-)\otimes\Q[y,y^{-1}] \to H_i(-)\otimes\Q[y,y^{-1},(1+y)^{-1}] =:H_*(-)\otimes\Q_{loc} \:.$$
This renormalization is needed to get for $M$ smooth the normalization condition
\begin{equation}
 T_{y*}(M)=MHT_{y*}(M)= T^*_y(TM)\cap [M] \in H_*(M)\otimes \Q[y]\:.
\end{equation}

It is the Hirzebruch class transformation $T_{y*}$, which {\em unifies} the (rationalized)
Chern class transformation  $c_*\otimes \Q$,  Todd class transformation $td_*$ and $L$-class transformation $L_*$
(compare \cite{BSY}). The corresponding characteristic number 
$$\chi_y(X):=deg\left(MHT_{y*}(X)\right) \in \Z[y]$$
for a singular (compact) algebraic variety $X$ captures information about the
{\em Hodge filtration} of Deligne's (\cite{De}) mixed Hodge structure on the rational cohomology (with compact support) $H_{(c)}^*(X;\Q)$ of $X$. In fact, by M.~Saito's work \cite{Sa1} one has an equivalence  
$$MHM(\{pt\}) \simeq \mhs^p$$
between mixed Hodge modules on a point space, and rational (graded) polarizable mixed Hodge structures.
Moreover, the corresponding mixed Hodge structure on rational cohomology with compact support
$$H_{c}^*(X;\Q)=H^*(\{pt\};k_!\Q^H_X)$$
(with $k: X\to \{pt\}$ a constant map) agrees with Deligne's one by another deep theorem of
M. Saito \cite{Sa2}. Therefore the transformations $MHC_y$ and $MHT_{y*}$ can be seen as a characteristic class version of the ring homomorphism
 $$\chi_y:K_0(\mhs^p) \to \Z[y,y^{-1}]$$ 
defined on the Grothendieck group of (graded) polarizable mixed Hodge structures by 
\begin{equation} \chi_y([H]):=\sum_p {\rm dim}\: Gr^p_F(H \otimes \C) \cdot (-y)^p \:,
\end{equation}
for $F$ the Hodge filtration of $H \in \mhs^p$.
Note that $\chi_y([\Lef])=-y$.\\

These characteristic class transformations are motivic refinements of the (rationalization of the) {\em Chern class transformation} $c_*\otimes \Q$ of MacPherson.
$MHT_{y*}$ factorizes by \cite{Sch3} as
$$MHT_{y*}: K_0(MHM(X)) \to   H_*(X)\otimes\Q[y,y^{-1}]\subset  H_*(X)\otimes \Q_{loc}\:,$$
fitting into a (functorial) commutative diagram
\begin{equation*}
\begin{CD}
F(X) @<\chi_{stalk}<< K_0(D^b_c(X)) @< rat << K_0(MHM(X))\\
@V c_*\otimes \Q VV  @V c_*\otimes \Q VV @VV MHT_{y*} V \\
H_*(X)\otimes \Q  @= H_*(X)\otimes \Q  @ < y=-1<<  H_*(X)\otimes \Q[y,y^{-1}] \:.
\end{CD}
\end{equation*}

Here $D^b_c(X)$ is the derived category of algebraically constructible sheaves  on $X$ (viewed as a complex analytic space), with
$rat$ associating to a (complex of) mixed Hodge module(s) the underlying perverse (constructible) sheaf complex,
 and $\chi_{stalk}$ is given by the Euler characteristic of the stalks.\\

Let us go back to the case when $X$ is a local complete intersection in some ambient smooth algebraic manifold. Then it is natural to compare  $cl_*(X)$ for a functorial
homology characteristic class theory $cl_*$ as above with the corresponding virtual characteristic class $cl_*^{\rm vir}(X)$. If $M$ is smooth,  then clearly we have that  $$cl_*^{\rm vir}(M)=cl^*(TM) \cap [M]=cl_*(M) \:.$$
However, if $X$ is singular, the difference between  the homology classes $cl_*^{\rm vir}(X)$ and $cl_*(X)$ depends in general on the singularities of $X$. This motivates the following
\begin{prob} 
Describe the difference $cl_*^{\rm vir}(X)-cl_*(X)$ in terms of the geometry of singular locus of $X$.
\end{prob}
The above problem is usually studied in order to understand the complicated homology classes $cl_*(X)$ in terms of the simpler virtual classes $cl_*^{\rm vir}(X)$, with the difference
terms measuring the complexity of the singularities of $X$.\\

This question was first studied for the Todd class transformation $td_*$,
where this difference term is vanishing. More precisely one has the
\begin{thm}[Verdier '76]\label{Verdier76} Assume that $i: X\to Y$ is regular embedding of (locally constant)
codimension $n$. Then the Todd class transformation $td_*$ commutes with specialization
(see \cite{V0}), i.e.
\begin{equation} i^!\circ td_*=  td_*\circ i^!: \;G_0(Y) \to H_{*-n}(X) \:. 
\end{equation}
\end{thm} 

Note that $Y$ need not be smooth.
\begin{cor}\label{cor-td}
Assume that $X$ can be realized as a local complete intersection in some ambient smooth algebraic manifold. Then $td_*^{\rm vir}(X)=td_*(X)$. Especially, if $X$ is  a global complete intersection given as the zero-fiber $X=\{f=0\}$ of a proper morphism $f: M\to \C^n$
on the algebraic manifold $M$, then the arithmetic genus
\begin{equation}
\chi(X)=\chi^{vir}(X)=\chi(X_t)
\end{equation}
of $X$ agrees with that of a nearby smooth fiber $X_t$ for $0 \neq t$ small and generic.
\end{cor}

The next case studied in the literature is the $L$-class transformation $L_*$ for $X$ a 
compact global complex hypersurface. 
\begin{thm}[Cappell-Shaneson '91] Assume $X$ is a global compact hypersurface $X=\{f=0\}$ for a proper complex algebraic function
$f: M\to \C$ on a complex algebraic manifold $M$. Fix a complex Whitney stratification of $X$ and let $\VV_0$ be the set of  strata $V$ with  ${\rm dim} V<{\rm dim} X$. Assume for simplicity, that all $V \in \VV_0$ are simply-connected (otherwise one has to use suitable twisted $L$-classes, see \cite{CS0, ShCa}). Then
\begin{equation}\label{L-formula}
L_*^{\rm vir}(X)-L_*(X)=\sum_{V \in \VV_0}\; \sigma({\rm lk}(V)) \cdot L_*(\bar{V}) \:,
\end{equation} 
where $\sigma({\rm lk}(V)) \in \Z$ is a certain signature invariant associated to the link pair of the stratum $V$ in $(M,X)$. 
\end{thm}

This result is in fact of topological nature, and holds more generally for a suitable compact stratified pseudomanifold $X$, which is PL-embedded into a manifold $M$ 
in real codimension two (see \cite{CS0, ShCa} for details). \\

If $cl_*=c_*$ is the Chern class transformation, the problem amounts to comparing the Fulton-Johnson class $c^{FJ}_*(X):=c_*^{vir}(X)$ (e.g., see \cite{Fu, FJ}) with the homology Chern class $c_*(X)$ of MacPherson. 
The difference between these two classes is measured by the so-called {\it Milnor class} $\MC_*(X)$ of $X$, which is studied in many references like 
\cite{Alu0, BLSS, BLSS2, BSS, Max, PP, Sch, Sch2, Y}. This is a homology class supported on the singular locus of $X$, and for a global hypersurface it was computed in \cite{PP} (see also \cite{Sch, Sch2, Y, Max}) as a weighted sum in the Chern-MacPherson classes of closures of singular strata of $X$, the weights depending only on the normal information to the strata. For example, if $X$ has only {\em isolated} singularities, the Milnor class equals (up to a sign) the sum of the local {\em Milnor numbers} attached to the singular points.
In the following section we explain our approach \cite{Sch, Sch2} through 
{\em nearby and vanishing cycles} (for constructible functions), which recently was adapted to the motivic Hirzebruch
and Chern class transformations \cite{CMSS, Sc}.

\section{Nearby and vanishing cycles}

Let us start to explain some basic constructions for constructible functions in the complex algebraic context
(compare \cite{SchB, Sch, Sch2}).
Here we work in the classical topology on the complex analytic space $X$
associated to a separated scheme of finite type over $Spec(\C)$.
\begin{defn}
A function $\alpha: X\to \Z$ is called (algebraically) constructible, if it satisfies one of
the following two equivalent properties:
\begin{enumerate}
\item $\alpha$ is a finite sum $\alpha = \sum_{j} n_{j}\cdot 1_{Z_{j}}$, with $n_{j}\in \Z$ and
$1_{Z_{j}}$ the characteristic function of the closed complex algebraic subset $Z_{j}$ of $X$.
\item   
$\alpha$ is (locally) constant on the strata of a
complex algebraic Whitney b-regular stratification of $X$.
\end{enumerate}
\end{defn}

This notion is closely related to the much more sophisticated notion of (algebraically) constructible (complexes 
of) sheaves on $X$. A sheaf $\FC$ of (rational) vector-spaces on $X$ with finite dimensional stalks
is (algebraically) constructible, if there
exists a complex algebraic Whitney b-regular stratification as above such that the restriction of $\FC$ to all
strata is locally constant. Similarly, a bounded complex of sheaves is constructible, if all it cohomology sheaves
have this property, and we denote by $D^{b}_{c}(X)$ the corresponding derived category of bounded constructible
complexes on $X$. The Grothendieck group of the triangulated category $D^{b}_{c}(X)$ is denoted by $K_0(D^{b}_{c}(X))$.\\

Since we assume that all stalks of a constructible complex are finite dimensional, by taking stalkwise
the Euler characteristic we get a natural group 
homomorphism 
\begin{equation}
\chi_{stalk}:  K_0(D^{b}_{c}(X)) \to F(X);\: [\FC]\mapsto \left(x\mapsto \chi(\FC_x)\right) \:.
\end{equation}
Here $F(X)$ is the group of (algebraically) constructible functions on $X$.
It is easy to show that natural transformation $\chi_{stalk}$ is surjective.\\

As is well known (and explained in detail in \cite{SchB}), all the usual functors in sheaf theory, which respect
the corresponding category of constructible complexes of sheaves, induce by the epimorphism $\chi_{stalk}$ well-defined
group homomorphisms on the level of constructible functions. We just recall these, which are important for later
applications or definitions. 
\begin{defn}\label{calc-F}
Let $f: X \to Y$ be an algebraic map of complex spaces associated
to separated schemes of finite type over $Spec(\C)$. Then one has the following transformations:\\
{\bf (1) pullback:} $f^{*}: F(Y)\to F(X);\; \alpha\mapsto \alpha\circ f$, which corresponds to the usual pullback
of sheaves
$$f^*: D^{b}_{c}(Y)\to D^{b}_{c}(X) \:.$$
{\bf (2) exterior product:}  $\alpha \boxtimes \beta \in F(X\times Y)$ for $\alpha\in F(X)$ and 
$\beta \in F(Y)$, given by $\alpha \boxtimes \beta((x,y)):=\alpha(x)\cdot \beta(y)$.
This corresponds on the sheaf level to the exterior product
$$\boxtimes^L: D^{b}_{c}(X) \times D^{b}_{c}(Y)\to D^{b}_{c}(X\times Y)\:.$$
{\bf (3) Euler characteristic:} Suppose $X$ is compact and $Y=\{pt\}$ is a point. Then one has 
$\chi: F(X)\to \Z$,
corresponding to 
$$R\Gamma(X,\cdot)=k_*: D^{b}_{c}(X)\to D^{b}_{c}(\{pt\}) $$ 
on the level of constructible complexes of sheaves, with $k: X\to \{pt\}$ the constant proper map. By linearity it is  
characterized by the convention that for a compact complex algebraic subspace $Z\subset X$ 
\begin{equation}
 \chi(1_{Z}) := \chi\left(H^{*}(Z;\Q)\right)
\end{equation}
is just the usual Euler characteristic of $Z$.\\
{\bf (4) proper pushdown:} Suppose $f$ is proper. Then one has $f_{*}=f_{!}: F(X)\to F(Y)$, corresponding to 
$$Rf_{*}=Rf_{!}:  D^{b}_{c}(X)\to D^{b}_{c}(Y)$$
on the level of constructible complexes of sheaves. Explicitly it is given by 
\begin{equation}
 f_{*}(\alpha)(y) := \chi(\alpha |_{\{f=y\}})\:,
\end{equation}
and in this form it goes back to the paper \cite{MP} of MacPherson.\\
{\bf (5) nearby cycles:} Assume $Y=\C$ and let $X_{0}:=\{f=0\}$ be the zero fiber. Then one has
$\psi_{f}: F(X)\to F(X_{0})$, corresponding to Deligne's nearby cycle functor
$$\psi_{f}: D^{b}_{c}(X)\to D^{b}_{c}(X_0) \:.$$
This was first introduced in \cite{V} by using resolution of singularities
(compare with \cite{SchB} for another approach using stratification theory).
By linearity, $\psi_{f}$ is uniquely defined by the convention that for a closed complex algebraic
subspace $Z\subset X$ the value 
\begin{equation}
 \psi_{f}(1_{Z})(x):=\chi\left( H^*(F_{f|_{Z},x};\Q)\right)
\end{equation}
is just the {\em Euler-characteristic of a local Milnor fiber} $F_{f|_{Z},x}$ of $f|_{Z}$ at $x$. 
Here this {\em local Milnor fiber} at $x$ is given by 
\begin{equation}
 F_{f|_{Z},x}:=Z\cap B_{\epsilon}(x)\cap \{f=y\} \:,
\end{equation}
with $0<|y|<<\epsilon<<1$ and $B_{\epsilon}(x)$ an open
(or closed) ball of radius $\epsilon$ around $x$ (in some local coordinates). 
Here we use the theory of a Milnor fibration of a function $f$ on the singular space $Z$ (compare \cite{Le, SchB}).\\
{\bf (6) vanishing cycles:} Assume $Y=\C$ and let $i: X_{0}:=\{f=0\}\hookrightarrow X$ be the inclusion
of the zero-fiber. Then one has
$\phi_{f}: F(X)\to F(X_{0});\;\phi_{f} := \psi_{f} - i^{*} $, corresponding to Deligne's vanishing cycle functor
$$\phi_{f}:  D^{b}_{c}(X)\to D^{b}_{c}(X_0) \:.$$
By linearity, $\phi_{f}$ is uniquely defined by the convention that for a closed complex algebraic
subspace $Z\subset X$ the value 
\begin{equation}
 \phi_{f}(1_{Z})(x):=\chi\left( H^*(F_{f|_{Z},x};\Q)\right)-1 = \chi\left( \tilde{H}^*(F_{f|_{Z},x};\Q)\right)
\end{equation}
is just the {\em reduced Euler-characteristic of a local Milnor fiber} $F_{f|_{Z},x}$ of $f|_{Z}$ at $x$.
\end{defn}

\begin{remark}
Let the global hypersurface $X=\{f=0\}$ be the zero-fiber of an algebraic function
$f: M\to \C$ on the complex algebraic manifold $M$. Then the support of $\phi_f(1_M)$
is contained in the {\em singular locus} $X_{sing}$ of $X$:
$$supp\left( \phi_f(1_M)\right) \subset X_{sing}\:.$$
And $\phi_f(1_M)|X_{sing}$ is (up to a sign) the {\em Behrend function} of $X_{sing}$ (see \cite{Beh}),
an intrinsic constructible function of the singular locus appearing in relation to
Donaldson-Thomas invariants.
\end{remark}

A beautiful result of Verdier \cite{V, Ken2} shows that for a global hypersurface MacPherson's Chern class transformation
$c_*$ commutes with {\em specialization}, if one uses the nearby cycle functor $\psi_f$ on the level of constructible functions
(as opposed to  the pullback functor $i^*$ for the corresponding inclusion $i: X=\{f=0\}\to Y$).

\begin{thm}[Verdier '81] \label{Verdier81} Assume that $X=\{f=0\}$ is a global hypersurface (of
codimension one) in $Y$ given by the zero-fiber of a complex algebraic function $f: Y \to \C$. 
Then the MacPherson Chern class transformation $c_*$ commutes with specialization
(see \cite{V, Ken}), i.e.
\begin{equation}
i^!\circ c_*=  c_*\circ \psi_f: F(Y) \to H_{*-1}(X) 
\end{equation}
for the closed inclusion $i: X=\{f=0\}\to Y$.
\end{thm} 
Note that $Y$ need not be smooth. As an immediate application one gets by (\ref{gysin-vir}) and  
(\ref{factor-N}) the following important result
(compare \cite{Sch, Sch2}):

\begin{cor}\label{chern-cor} Assume that $X=\{f=0\}$ is a global hypersurface (of
codimension one) in some ambient smooth algebraic manifold $M$,
given by the zero-fiber of a complex algebraic function $f: M \to \C$.
Then
\begin{equation}
 c^{vir}_*(X)-c_*(X)=c_*(\psi_f(1_M))-c_*(1_X)=c_*(\phi_f(1_M))\in H_*(X_{sing}) \:,
\end{equation}
since $supp(\phi_f(1_M))\subset X_{sing}$. Here we also use the naturality of $c_*$ for the closed
inclusion $X_{sing}\to X$ to view this difference term as a {\em localized class} in $H_*(X_{sing})$.
In particular:
\begin{enumerate}
 \item $c^{vir}_i(X)=c_i(X)\in H_i(X)$ for all $i>dim\;X_{sing}$.
\item  If $X$ has only isolated singularities (i.e. $dim\;X_{sing}=0$), then
$$c^{vir}_*(X)-c_*(X) =\sum_{x \in X_{sing}} \;\chi\left(\tilde{H}^*(F_x;\Q)\right)\:,$$
where  $F_x$ is the local Milnor fiber of the isolated hypersurface singularity $(X,x)$.
\item If $f: M\to \C$ is proper, then
$$deg\left(c_*(\phi_f(1_M))\right) = deg\left( c^{vir}_*(X)-c_*(X)\right)=
\chi(X_t)-\chi(X)$$
is the difference between the Euler characteristic of a global nearby smooth fiber $X_t=\{f=t\}$
(for $0\neq |t|$ small enough) and of the special fiber $X=\{f=0\}$.
\end{enumerate}
\end{cor}

For a general local complete intersection $X$ in some ambient smooth algebraic manifold 
(e.g. a local hypersurface of codimension one), one doesn't have global equations so that
the theory of nearby and vanishing cycles can't be applied directly. Instead one has to combine
them with the {\em deformation to the normal cone} leading to Verdier's theory of
{\em specialization functors} (compare \cite{Sch, Sch2}). But even if $X=\{f=0\}$ is a global
complete intersection inside the ambient smooth algebraic manifold $M$, given by the zero-fiber of a
complex algebraic map $f: M\to \C^n$, one doesn't have a theory of nearby and vanishing cycles, because
a local theory of Milnor fibers for $f$ is missing (if $n>1$). 
But if one fixes an {\em ordering} of the components of $f$ (or of the coordinates on $\C^n$),
then a corresponding {\em local Milnor fibration} exists for any ordered tuple
$$(f):=(f_1,\dots,f_n): Z\to \C^n$$
of complex algebraic functions on the singular algebraic variety $Z$ (as observed in \cite{MCP}).

\begin{defn}[Nearby and vanishing cycles for an ordered tuple]
Let $(f):=(f_1,\dots,f_n): Y \to \C^n$ be an ordered $n$-tuple of complex algebraic functions on $Y$,
with $X:=\{f=0\}=\{f_1=0,\dots,f_n=0\}$ the zero-fiber of $(f)$.
Then {\em nearby cycles} of $(f):=(f_1,\dots,f_n)$ are defined by iteration as 
\begin{equation}
 \psi_{(f)}:= \psi_{f_1}\circ \cdots \circ \psi_{f_n}:  F(Y)\to F(X)\:.
\end{equation}
By linearity, $\psi_{(f)}$ is uniquely defined by the convention that for a closed complex algebraic
subspace $Z\subset Y$ the value 
\begin{equation}
 \psi_{(f)}(1_{Z})(x):=\chi\left( H^*(F_{(f)|_{Z},x};\Q)\right)
\end{equation}
is just the {\em Euler-characteristic of a local Milnor fiber} $F_{(f)|_{Z},x}$ of $(f)|_{Z}$ at $x$. 
Here this {\em local Milnor fiber of $(f)$} at $x$ is given by 
\begin{equation}
 F_{(f)|_{Z},x}:=Z\cap B_{\epsilon}(x)\cap \{f_1=y_1,\dots ,f_n=y_n\} \:,
\end{equation}
with $0<|y_n|<<\cdots << |y_1| <<\epsilon<<1$ and $B_{\epsilon}(x)$ an open
(or closed) ball of radius $\epsilon$ around $x$ (in some local coordinates, compare \cite{MCP}).

The corresponding {\em vanishing cycles} of $(f)$ are defined by
\begin{equation}
\phi_{(f)} := \psi_{(f)} - i^{*}:  F(Y)\to F(X)\:,
\end{equation}
with $i: X\to Y$ the closed inclusion.
By linearity, $\phi_{(f)}$ is uniquely defined by the convention that for a closed complex algebraic
subspace $Z\subset X$ the value 
\begin{equation}
 \phi_{(f)}(1_{Z})(x):=\chi\left( H^*(F_{(f)|_{Z},x};\Q)\right)-1 = \chi\left( \tilde{H}^*(F_{(f)|_{Z},x};\Q)\right)
\end{equation}
is just the {\em reduced Euler-characteristic of a local Milnor fiber} $F_{(f)|_{Z},x}$ of $(f)|_{Z}$ at $x$.
\end{defn}

Note that again $supp\left(\phi_{(f)}(1_{M})\right) \subset X_{sing}$ in case the ambient space $Y=M$ is a smooth algebraic
manifold. Assume moreover that $X$ is of codimension $n$ so that the regular embedding $i: X\to Y$ factorizes into $n$
regular embeddings of codimension one $i=i_n\circ \cdots \circ i_1$:
\begin{equation}\begin{CD}
 X= \{f_1=0,\dots,f_n=0\} @>i_1>> \{f_2=0,\dots,f_n=0\}  @>i_2>> \cdots\\
\cdots   \{f_{n-1}=0,f_n=0\}   @> i_{n-1} >>  \{f_n=0\}@> i_n >> Y \:.
\end{CD}\end{equation}
By the functoriality of the Gysin homomorphisms one gets
$$i^!=i_1^!\circ \cdots \circ i_n^!: H_*(Y)\to H_{*-n}(X)\:.$$
Since in Verdier's spezialisation theorem (\ref{Verdier81}) the ambient space need not be smooth,
we can apply it inductively to all embeddings $i_j$ (for $j=n,\dots,1$) above.

\begin{cor}\label{ordered-cor} Assume that $X=\{f=0\}=\{f_1=0,\dots,f_n=0\}$ is a global complete intersection (of
codimension $n$) in some ambient smooth algebraic manifold $M$,
given by the zero-fiber of an ordered $n$-tuple of complex algebraic function $(f):=(f_1,\dots,f_n): M \to \C^n$.
Then
\begin{equation}
 c^{vir}_*(X)-c_*(X)=c_*(\phi_{(f)}(1_M))\in H_*(X_{sing}) \:,
\end{equation}
since $supp(\phi_{(f)}(1_M))\subset X_{sing}$. Here we also use the naturality of $c_*$ for the closed
inclusion $X_{sing}\to X$ to view this difference term as a {\em localized class} in $H_*(X_{sing})$.
In particular:
\begin{enumerate}
 \item $c^{vir}_i(X)=c_i(X)\in H_i(X)$ for all $i>dim\;X_{sing}$.
\item  If $X$ has only isolated singularities (i.e. $dim\;X_{sing}=0$), then
$$c^{vir}_*(X)-c_*(X) =\sum_{x \in X_{sing}} \;\chi\left(\tilde{H}^*(F_x;\Q)\right)\:,$$
where  $F_x$ is the local Milnor fiber of the ordered $n$-tuple $(f)$ at the isolated singularity $x$.
\item If $(f)=(f_1,\dots,f_n): M\to \C^n$ is proper, then
$$deg\left(c_*(\phi_{(f)}(1_M))\right) = deg\left( c^{vir}_*(X)-c_*(X)\right)=
\chi(X_t)-\chi(X)$$
is the difference between the Euler characteristic of a global nearby smooth fiber $X_t=\{f_1=t_1,\dots, f_n=t_n\}$
(for $t=(t_1,\dots ,t_n)$ with $0<|t_n|<<\cdots << |t_1|$
 small enough) and of the special fiber $X=\{f=0\}$.
\end{enumerate}
\end{cor}

As explained in section $2$, the motivic Hirzebruch and Chern class transformations
$T_{y*}, MHT_{y*}$ and $mC_y, MHC_y$ can be seen as ``motivic or Hodge theoretical liftings'' of
the (rationalized) Chern class transformation $c_*$ under the comparison maps
$$\begin{CD}
K_0(var/Y) @> \chi_{Hdg} >> K_0(MHM(Y)) @> rat >> K_0(D^b_c(Y)) @> \chi_{stalk} >> F(Y) \:.
  \end{CD}$$
Here these Grothendieck groups have the same calculus as for constructible functions in
definition \ref{calc-F}(1-4), respected by these comparison maps.
So it is natural to try 
to extend known results about MacPherson's Chern class transformation $c_*$ to these transformations.
In the ``motivic'' (resp. `` Hodge theoretical'') context this has been worked out in \cite{BSY}
(resp. \cite{Sch3}) for
\begin{enumerate}
 \item the functorialty under {\em push down} for proper algebraic morphism.
\item the functorialty under {\em exterior products}.
\item the functorialty under {\em smooth pullback} given by a related {\em Verdier-Riemann-Roch theorem}.
\end{enumerate}

And recently we could also prove the ``counterpart'' of Verdier's specialization theorem (\ref{Verdier81}).
Let $X=\{f=0\}$ be a global hypersurface in $Y$ given by the zero-fiber of a complex algebraic function $f$ on $Y$:
$$\begin{CD}
X:=\{f=0\} @> i >> Y @> f >> \C \:.
\end{CD}$$
First note that one can use the nearby and vanishing cycle functors $\psi_f$ and $\phi_f$ either on the motivic level of 
localized relative Grothendieck groups $$\MC(var/-):=K_0(var/-)[\Lef^{-1}]$$ (see \cite{Bit, GLM}), 
or on the Hodge-theoretical level of algebraic mixed Hodge modules 
(\cite{Sa0, Sa1}), ``lifting'' the corresponding functors on the level of algebraically constructible sheaves 
(\cite{Sch}) and algebraically constructible functions as introduced before, so that the following diagram commutes:

\begin{equation}\label{all-nearby}\begin{CD}
\MC(var/Y) @> \psi_f^m, \phi_f^m >>  \MC(var/X)\\
@V \chi_{Hdg} VV  @V \chi_{Hdg} VV  \\
K_0(\mh(Y)) @> \psi'^H_f, \phi'^H_f >>  K_0(\mh(X))\\
@V rat VV @V rat VV \\
K_0(D^b_c(Y)) @> \psi_f, \phi_f >>  K_0(D^b_c(X)) \\
@V \chi_{stalk} VV @V \chi_{stalk} VV \\
F(Y) @> \psi_f, \phi_f >> F(X)  \:.
\end{CD}\end{equation}
We also use the notation $\psi'^H_f:= \psi^H_f[1]$ and $\phi'^H_f:=\phi^H_f[1]$ for the shifted
functors, with $\psi^H_f, \phi^H_f: MHM(Y)\to MHM(X)$ and $\psi_f[-1],\phi_f[-1]: Perv(Y)\to Perv(X)$ preserving mixed Hodge modules and perverse sheaves, respectively. 
On the level of Grothendieck groups one simply has $\phi_f^m =\psi_f^m-i^*$ and $ \phi'^H_f=\psi'^H_f-i^*$.

\begin{remark}\rm
The motivic nearby and vanishing cycles functors of \cite{Bit, GLM}
take values in a refined {\em equivariant} localized Grothendieck group
$\MC^{\hat{\mu}}(var/X)$ of equivariant algebraic varieties over $X$ with a ``good'' action
of the pro-finite group $\hat{\mu}=\lim \mu_n$ of  roots of unity.
For mixed Hodge modules this corresponds to an action of the semi-simple part of the monodromy.
But in the following applications we don't need to take this action into account. 
Also note that for the commutativity of diagram (\ref{all-nearby}) one has to use $\psi'^H_f, \phi'^H_f$ (as opposed to 
$\psi^H_f, \phi^H_f$).
\end{remark}

Now we are ready to formulate the main new result from \cite{Sc}.

\begin{thm}[Sch\"{u}rmann '09]\label{Sch09}
Assume that $X=\{f=0\}$ is a global hypersurface of codimension one given by the zero-fiber of a complex algebraic
function $f: Y\to \C$.
Then the motivic Hodge-Chern class transformation $MHC_y$ commutes with specialization in the following sense:
\begin{equation}
(1+y)\cdot MHC_y(\; \psi'^H_f (-)\;) =  i^!MHC_y(-)   
\end{equation} 
as transformations $K_0(MHM(Y))\to G_0(X)[y,y^{-1}]$.
\end{thm}

Again the smoothness of $Y$ is not needed. The appearence of the factor $(1+y)$ should not be a surprise,
as it can already be seen in the case of a smooth hypersurface $X$ inside a smooth ambient manifold $Y$,
$$(1+y)\cdot MHC_y(X)= i^!MHC_y(Y) \:,$$
if one recalls (\ref{gysin-vir}),  (\ref{factor-N}) and 
the normalization condition (\ref{C-norm}), with 
$$\Q^H_X = i^*\Q^H_Y\simeq \psi_f'^H(\Q^H_Y)$$ 
in this special case.
But the proof of this theorem given in \cite{Sc} is far away from the geometric applications
described here. In fact it uses the algebraic theory of nearby and vanishing cycles in the context of $D$-modules
given by the {\em $V$-filtration of Malgrange-Kashiwara}, together with a specialization result about the {\em filtered de Rham complex}
of the filtered $D$-module underlying a mixed Hodge module.\\
 
Using Verdier's result that the Todd class transformation $td_*$ commutes with specialization (see theorem \ref{Verdier76}), one gets (\cite{Sc}):
\begin{cor} \label{main-sch}
Assume that $X=\{f=0\}$ is a global hypersurface of codimension one given by the zero-fiber of a complex algebraic
function $f: Y\to \C$. Then the motivic Hirzebruch class transformation
${MHT_y}_*$ commutes with specialization, that is: 
\begin{equation}\label{sp}
{MHT_y}_*(\psi'^H_f(-))=i^! {MHT_y}_*(-) 
\end{equation}
as transformations $K_0({\rm MHM}(Y)) \to H_*(X) \otimes \Q[y,y^{-1}]$.
\end{cor}
Again the smoothness of $Y$ is not needed here, but only the fact that $X=\{f=0\}$ is a global hypersurface (of codimension one) is needed.
Also the factor $(1+y)$ 
in theorem \ref{Sch09} cancelled out by the renormalization factor $(1+y)^{-i}\cdot $ on $H_i(-)$ used in the definition of
$MHT_{y*}$, since the Gysin map $i^!: H_*(Y)\to H_{*-1}(X)$ shifts this degree by one.\\

By the definition of $\psi_f^m$ in \cite{Bit,GLM} one has that
$$\psi_f^m(K_0(var/Y))\subset im(K_0(var/X)\to \MC(X))\:,$$
so ${MHT_{y}}_*\circ \psi_f^m$ maps $K_0(var/Y)$ into
$H_*(X)\otimes\Q[y]\subset H_*(X)\otimes\Q[y,y^{-1}]$.
Together with \cite{Sch3}[Prop.5.2.1] one therefore gets the following commutative diagram:

\begin{equation}\label{special}\begin{CD}
K_0(var/Y) @> {T_{y}}_*\circ \psi_f^m = > i^!\circ {T_{y}}_* > H_*(X)\otimes \Q[y]\\
@V \chi_{Hdg} VV @VVV\\
K_0(\mh(Y)) @> {MHT_{y}}_*\circ \psi'^H_f = > i^!\circ {MHT_{y}}_* > H_*(X)\otimes \Q[y,y^{-1}]\\
@V \chi_{stalk}\circ rat VV @VV y=-1 V \\
F(Y) @> c_*\circ \psi_f = > i^! c_* > H_*(X)\otimes \Q \:.
\end{CD}\end{equation}

As before one gets the following result from Theorem \ref{Sch09} and Corollary \ref{main-sch} together with (\ref{gysin-vir}) and  
(\ref{factor-N}):

\begin{lemma}\label{L1} Assume that $X=\{f=0\}$ is a global hypersurface (of
codimension one) in some ambient smooth algebraic manifold $M$,
given by the zero-fiber of a complex algebraic function $f: M \to \C$.
Then
\begin{equation} mC_y^{\rm vir}(X)= mC_{y}\left(\psi_f^m(\left[id_M\right])\right)=
MHC_{y*}\left(\psi'^H_f(\left[\Q^H_M\right])\right) \:,
\end{equation}
and
\begin{equation} T_{y*}^{\rm vir}(X)= T_{y*}\left(\psi_f^m(\left[id_M\right])\right)=
MHT_{y*}\left(\psi'^H_f(\left[\Q^H_M\right])\right) \:.
\end{equation}
\end{lemma}

If $i: X=\{f=0\}\to M$ is the closed inclusion, then one has $i^*([id_M])=[id_X]$ and $i^*([\Q^H_M])=[\Q^H_X]$.
So by $\phi_f^m =\psi_f^m-i^*$ and $ \phi'^H_f=\psi'^H_f-i^*$ (on the level of Grothendieck groups) one gets
(compare \cite{CMSS}):

\begin{cor}\label{main-cor} Assume that $X=\{f=0\}$ is a global hypersurface (of
codimension one) in some ambient smooth algebraic manifold $M$,
given by the zero-fiber of a complex algebraic function $f: M \to \C$.
Then
\begin{equation}\begin{split}
 mC_y^{vir}(X)-mC_y(X) &= mC_{y}\left(\phi_f^m(\left[id_M\right])\right) \\
&= MHC_{y*}\left(\phi'^H_f(\left[\Q^H_M\right])\right)\in G_0(X_{sing})[y]\:,
\end{split}\end{equation}
and
\begin{equation} \label{phi-Ty} \begin{split}
T_{y*}^{\rm vir}(X) - T_{y*}(X) &=T_{y*}\left(\phi_f^m(\left[id_M\right])\right) \\
&= MHT_{y*}\left(\phi'^H_f(\left[\Q^H_M\right])\right) \in H_*(X_{sing})\otimes \Q[y]\:.
\end{split}\end{equation}
Here we use  
$$supp\left( \phi'^H_f\left(\Q^H_M\right) \right)\subset X_{sing}$$ 
and the naturality of our characteristic class transformations for the closed
inclusion $X_{sing}\to X$.
In particular:\\
$(1)$ $T_{y,i}^{vir}(X)=T_{y,i}(X)\in H_i(X)\otimes \Q[y]$ for all $i>dim\;X_{sing}$.\\
$(2)$  If $X$ has only isolated singularities (i.e. $dim\;X_{sing}=0$), then
\begin{equation} \begin{split}
mC^{vir}_y(X)-mC_y(X) &=\sum_{x \in X_{sing}} \;\chi_y\left(\tilde{H}^*(F_x;\Q)\right) \\
&= T_{y*}^{\rm vir}(X) - T_{y*}(X)\:,
\end{split}\end{equation}
where  $F_x$ is the Milnor fiber of the isolated hypersurface singularity $(X,x)$.\\
$(3)$ If $f: M\to \C$ is proper, then
\begin{equation}\begin{split}
deg\left( MHC_{y*}\left(\phi'^H_f(\left[\Q^H_M\right])\right) \right) &= 
\chi_y\left(H^*(X_t;\Q)\right)-\chi_y\left(H^*(X;\Q)\right) \\
&= deg\left( MHT_{y*}\left(\phi'^H_f(\left[\Q^H_M\right])\right)  \right)
\end{split}\end{equation}
is the difference between the $\chi_y$-characteristics of a global nearby smooth fiber $X_t=\{f=t\}$
(for $0\neq |t|$ small enough) and of the special fiber $X=\{f=0\}$.
\end{cor}

\begin{remark}\label{spec} (Hodge polynomials vs. Hodge spectrum) \rm
Let us explain the precise relationship between the Hodge spectrum and the less-studied $\chi_y$-polynomial of the Milnor fiber of a hypersurface singularity. Here 
we follow notations and sign conventions similar to those in \cite{GLM}. Denote by $\mhs^{\rm mon}$ the abelian category of mixed Hodge structures endowed with an automorphism of finite order, and by $K^{mon}_0(\mhs)$ the corresponding Grothendieck ring. There is a natural linear map called the {\it Hodge spectrum}, 
$${\rm hsp}: K^{\rm mon}_0(\mhs) \to \Z[\Q] \simeq \bigcup_{n \geq 1} \Z[t^{1/n},t^{-1/n}] \:,$$ 
such that
\begin{equation}\label{hsp}
{\rm hsp}([H]):=\sum_{ \alpha \in \Q \cap [0,1)} t^{\alpha} \left( \sum_{p \in \Z} {\rm dim}(Gr^{p}_F H_{\C,\alpha}) t^p \right)
\end{equation}
for any mixed Hodge structure $H$ with an automorphism $T$ of finite order, where $H_{\C}$ is the underlying complex vector space of $H$, $H_{\C,\alpha}$ is the eigenspace of $T$ with eigenvalue 
$\exp(2\pi i \alpha)$, and $F$ is the Hodge filtration on $H_{\C}$.
It is now easy to see that the $\chi_y$-polynomial ${\rm \chi_y}([H])$ of $H$ is obtained from ${\rm hsp}([H])$ by substituting $t=1$ in $t^{\alpha}$ for $\alpha \in \Q \cap [0,1)$ and  $t=-y$ in $t^p$ for $p\in \Z$.
\end{remark}

As already explained before, Corollary \ref{main-cor} reduces for the value $y=-1$ of the parameter to the (rationalized version of) Corollary \ref{chern-cor}.
Since the ambient space in Theorem \ref{Sch09} and Corollary \ref{main-sch} need not be smooth, one can generalize in the same way the Corollary \ref{ordered-cor}
for a global complete intersection $X=\{f=0\}=\{f_1=0,\dots,f_n=0\}$ (of
codimension $n$) in some ambient smooth algebraic manifold $M$,
given by the zero-fiber of an ordered $n$-tuple of complex algebraic function $(f):=(f_1,\dots,f_n): M \to \C^n$. Here we leave the details to the reader.\\

It is also very interesting to look at the other specializations of Corollary \ref{main-cor}  for $y=0$ and $y=1$.
Let us first consider the case when $y=0$. Note that in general $T_{0*}(X)\neq td_*(X)$ for a singular complex algebraic variety (see \cite{BSY}).
But if $X$ has only {\em Du Bois singularities} (e.g., rational singularities, cf. \cite{Sa2}), then by \cite{BSY} we have
$T_{0*}(X)= td_*(X)$. 
So if a global hypersurface $X=\{f=0\}$ has only {\em Du Bois singularities}, then by Corollaries \ref{cor-td} and \ref{main-cor} we get:
$$MHT_{0*}\left(\phi'^H_f(\left[\Q^H_M\right])\right)=0 \in H_*(X)\otimes \Q \:.$$ 

This vanishing (which is in fact a class version of Steenbrink's {\em cohomological insignificance of $X$} \cite{St}) imposes interesting geometric identities on the corresponding Todd-type invariants of the singular locus. For example, we obtain the following
\begin{cor}  If the global hypersurface $X$ has only {\em isolated Du Bois singularities}, then
\begin{equation}\label{van} {\rm dim}_{\C} Gr^0_FH^n(F_x;\C)=0\end{equation}
for all  $x\in X_{sing}$, with $n=dim\;X$. 
\end{cor}
\noindent It should be pointed out that in this setting a result of Ishii \cite{I} implies that (\ref{van}) is in fact equivalent to $x\in X_{sing}$ being an isolated Du Bois hypersurface singularity.
Also note that in the arbitrary singularity case, the {\em Milnor-Todd class } 
$$ T_{0*}\left(\phi_f^m(\left[id_M\right])\right) 
= MHT_{0*}\left(\phi'^H_f(\left[\Q^H_M\right])\right) \in H_*(X_{sing})\otimes \Q$$
carries interesting non-trivial information about the singularities of the hypersurface $X$.\\

Finally, if $y=1$, the formula (\ref{phi-Ty}) should be compared to the Cappell-Shaneson topological result of (\ref{L-formula}).
While it can be shown (compare with \cite{Max05}) that the normal contribution $\sigma({\rm lk}(V))$ in (\ref{L-formula}) for a singular stratum $V \in \VV_0$ is in fact the 
{\em signature} $\sigma(F_v)$ ($v \in V$) of the Milnor fiber (as a manifold with boundary) of the singularity in a transversal slice to $V$ in $v$, the precise relation
between $\sigma(F_v)$ and $\chi_1(F_v)$ is in general very difficult to understand.
For $X$ a {\em rational homology manifold}, one would like to have
a  ``local Hodge index formula'' 
$$\sigma(F_v)\stackrel{?}{=}\chi_1(F_v) \:,$$
which is presently not available. But if the hypersurface $X$ is a {\em rational homology manifold with only isolated
singularities}, then this expected equality follows from \cite{St5}[Thm.11]. One therefore gets in this case (by a comparison
of the different specialization results for $L_*$ and $T_{1*}$)
the following conjectural interpretation of $L$-classes from \cite{BSY}
(see \cite{CMSS} for more details):
\begin{thm} Let $X$ be a compact complex algebraic variety with only isolated singularities,
which is moreover a  rational homology manifold and can be realized as a global hypersurface (of codimension one) in a complex algebraic manifold. Then
\begin{equation}
L_*(X)=T_{1*}(X)\in H_{2*}(X;\Q)\:.
\end{equation}
\end{thm}

\end{document}